# ON DIFFEOLOGY OF ORBIT SPACES

SERAP GÜRER AND PATRICK IGLESIAS-ZEMMOUR

*Dedicated to Richard S. Palais on the occasion of his 94th birthday.*
*This work is presented in deep appreciation of his foundational contributions*
*to the theory of transformation groups, which form the essential framework of this paper.*

ABSTRACT. We investigate the correspondence between the geometry of a smooth compact Lie group action on a manifold M and the intrinsic smooth structure of the orbit space M/G. While the action on M is classically organized by the orbit-type stratification, we show this structure fails to predict the intrinsic Klein stratification of the quotient, which partitions the space into the orbits of local diffeomorphisms, thereby classifying the space by its intrinsic singularity types. The correct correspondence, we prove, is governed by a finer partition on M: the *isostabilizer decomposition*. We establish a surjective map from the components of this partition to the Klein strata of M/G. As a corollary, we obtain by pullback a new canonical stratification on M, the *Inverse Klein Stratification*, and clarify its relationship with classical structures.

## INTRODUCTION

The study of spaces arising from the smooth action of a compact Lie group G on a manifold M is a cornerstone of differential geometry, with a rich history built upon the work of Palais, Bredon, and Mostov. A central challenge is that the orbit space M/G is generally a singular space. Classical approaches focused on imposing structure on the source manifold M, leading to the celebrated *orbit-type stratification* [Bre72, Pal60]. This partition has excellent properties—its strata are G-invariant submanifolds and it satisfies the frontier condition, making it a true stratification of M. For decades, it has been the primary tool for organizing the structure of a G-manifold. The projection of this partition to the orbit space M/G yields a true stratification of the quotient. However, this projected structure is extrinsic; it is inherited from the ambient manifold M rather than being defined by the intrinsic geometry of the quotient space itself. It answers the question of where points came from, but not necessarily what they "look like" locally.

The framework of diffeology offers a paradigm shift. It provides a powerful, intrinsic theory of smoothness for the quotient space M/G itself, where the logic of the category

*Date*: July 23, 2025.
2020 *Mathematics Subject Classification.* 57S15, 58A40, 58A35.
*Key words and phrases.* Diffeology, Orbit Space, Isostabilizer Decomposition, Klein Stratification, Orthofold.
(P.I-Z thanks the Hebrew University of Jerusalem, Israel, for its continuous academic support. He is also grateful for the stimulating discussions and assistance provided by the AI assistant Gemini (Google).)

takes precedence. This intrinsic geometry is captured by the *Klein stratification*, the partition of M/G into its local diffeomorphism orbits, which are the true measure of its singularities. This raises a fundamental question: what partition of M maps coherently to the Klein stratification of its quotient?

One might hope the classical orbit-type stratification would be the answer, but our key examples show this is not the case. The correspondence fails in subtle and profound ways, revealing that the very properties that make the orbit-type partition a good stratification of M make it the wrong tool for understanding M/G.[1]

We show in this paper that the correct tool is the finer *isostabilizer decomposition*, whose elements are the connected components of submanifolds where the stabilizer subgroup is constant.

We establish the main result of this paper: the *Correspondence Theorem*, which states that the canonical projection induces a surjective map from the isostabilizer decomposition of M to the Klein stratification of M/G. This provides the definitive dictionary between the local data of the action—the stabilizer and its slice representation—and the intrinsic singularity types of the quotient.

As a consequence of our Correspondence Theorem, we are able to define a new, canonical stratification on M itself—the *Inverse Klein Stratification*—by pulling back the Klein strata from the quotient. This allows us to clarify the full hierarchy of partitions on M: we show that the isostabilizer decomposition is the finest of all, refining both the classical orbit-type stratification and our new Inverse Klein Stratification, and that these latter two are, in general, incomparable. Our analysis also highlights the profound feature that the intrinsic smooth identity of M/G, which we term an *orthofold*, can be independent of the specific global action that created it.

The paper is organized as follows: Section I provides preliminaries on diffeology. Section II discusses the diffeological structure of M/G. Section III presents the illustrative examples that motivate our analysis. Section IV defines the isostabilizer decomposition. Section V establishes our main theorem and, as a corollary, constructs the Inverse Klein Stratification. Finally, Section VI provides concluding remarks.

I. PRELIMINARIES ON DIFFEOLOGY

In this section, we provide a brief introduction to the fundamental concepts of diffeology necessary for understanding the structure of orbit spaces as diffeological spaces. For a comprehensive treatment, we refer the reader to [PIZ13].

**1. Diffeological Spaces.**

---

[1] This journey from a natural but flawed conjecture to a more refined theorem, guided by the analysis of counterexamples, is a real-world instance of the process of "proofs and refutations" described by Imre Lakatos. Our initial difficulties, and indeed the conventional wisdom of focusing on the orbit-type stratification, stemmed from the "prejudice" that a well-behaved structure on the source manifold ought to be the primary analytical tool.

A *parametrization* in a set X is a map $P : U \to X$, where U is an open subset of some Euclidean space $\mathbf{R}^k$. A *diffeology* on a set X is a collection of parametrizations, called the *plots of the diffeology*, satisfying three axioms:

(1) **Covering:** Any constant map $P : U \to X$ is a plot.
(2) **Locality:** If $P : U \to X$ is a parametrization such that for every $r \in U$, there exists an open neighborhood $V \subseteq U$ of $r$ such that the restriction $P|_V$ is a plot, then P is a plot.
(3) **Smooth Compatibility:** If $P : U \to X$ is a plot and $f : V \to U$ is a smooth map between open subsets of Euclidean spaces, then the composition $P \circ f : V \to X$ is a plot.

A *diffeological space* is a set X equipped with a diffeology. Euclidean domains and smooth manifolds are first examples of diffeological spaces, where the plots are precisely the smooth maps (in the usual sense) from open subsets of Euclidean spaces into Euclidean spaces or manifold.

## 2. D-Topology.

Every diffeological space X carries a natural topology, called the *D-topology*. A subset $\mathcal{O} \subseteq X$ is defined to be *D-open* if for every plot $P : U \to X$, the preimage $P^{-1}(\mathcal{O})$ is an open subset of U in the standard Euclidean topology.

## 3. Subspaces and Quotients.

Given a diffeological space X and a subset $A \subseteq X$, the *subset diffeology* on A is defined by taking as plots all parametrizations $P : U \to A$ such that P, viewed as a parametrization into X, is a plot of X.

Given a diffeological space X and an equivalence relation $\sim$ on X, the *quotient diffeology* on the set of equivalence classes $X/\!\!\sim$ is defined as follows: a parametrization $\bar{P} : U \to X/\!\!\sim$ is a plot if for every $r \in U$, there exists an open neighborhood $V \subseteq U$ of $r$ and a plot $P : V \to X$ such that $\bar{P}|_V$ is the composition of P with the canonical projection $\pi : X \to X/\!\!\sim$.

The canonical projection $\pi : X \to X/\!\!\sim$ is a fundamental example of a *subduction*. A map $f : X \to Y$ between diffeological spaces is a *subduction* if it is surjective and a parametrization $P : U \to Y$ is a plot of Y if and only if for every $r \in U$, there exists an open neighborhood $V \subseteq U$ of $r$ and a plot $\tilde{P} : V \to X$ such that $P|_V = f \circ \tilde{P}$. Subductions are the epimorphisms in the category of diffeological spaces.

## 4. Local Diffeomorphisms and Invariants.

A map $f : X \to Y$ between diffeological spaces is *smooth* if for every plot P of X, the composition $f \circ P$ is a plot of Y. A smooth map $f$ is a *diffeomorphism* if it is bijective and its inverse $f^{-1}$ is also smooth.

A *local diffeomorphism* from X to Y is a map $f$ between two D-open subsets $U \subseteq X$ and $V \subseteq Y$ which is a diffeomorphism from U (with the subspace diffeology) to V (with the subspace diffeology). The collection of all local diffeomorphisms of a space X forms a pseudogroup, denoted $\text{Diff}_{\text{loc}}(X)$.

The diffeology of a space X determines several intrinsic invariants, in particular:

**Definition (Dimension Map).** *The* diffeological dimension *at a point $x \in X$, denoted* $\dim_x(X)$, *is an intrinsic invariant that measures the minimal dimension of a Euclidean space from which one can locally generate all plots around $x$. It is defined precisely as a minimax over families of plots pointed at $x$.*

This definition was introduced in [PIZ07] and is further elaborated with many examples in [PIZ13]. It coincides with the standard dimension for manifolds.

**Definition (Stratification).** *A* stratification *of a diffeological space* X *is a partition into disjoint subsets (strata) satisfying the* frontier condition *with respect to the D-topology. This is a* basic *condition:*

   [B] *The closure of every stratum* $S \in \mathscr{S}$, *in the D-topology, is a union of strata.*

Diffeological stratifications may satisfy other specific conditions, often represented by various labels (e.g., [LF] for locally finite, [M] for manifold strata). For a presentation of the general theory of stratifications in diffeology, we refer to [GIZ23].

In the spirit of Klein's Erlangen Program, the local geometry of a diffeological space is classified by its pseudogroup of local diffeomorphisms, $\mathrm{Diff}_{\mathrm{loc}}(X)$. The orbits of this pseudogroup partition the space X into sets of points that are locally diffeomorphic to each other. These orbits are called the *Klein strata* of X [PIZ13, §1.42]. They form the finest partition of X that respects its local smooth structure.

**Proposition (Klein Stratification is Basic).** *The partition of a diffeological space into Klein strata is a basic stratification: The closure of a Klein stratum is a union of Klein strata.*

*Proof.* This proposition is a standard result in diffeology (see [PIZ13, §1.42] and [PIZ25, Section 74]). The proof relies on the fact that local diffeomorphisms are also homeomorphisms with respect to the D-topology. Since the Klein strata are the orbits of this pseudogroup, the closure of a stratum is necessarily a union of strata, thus satisfying the frontier condition. □

Points in the same Klein stratum share the same local geometry, understood as local diffeological structure. We can also say that the Klein stratification partitions the space into points with the same types of singularities. In particular,

**Proposition (Klein Stratification and Dimension).** *The diffeological dimension is constant on the Klein strata.*

*Proof.* By definition, any two points in the same Klein stratum are locally diffeomorphic. Since the diffeological dimension is an invariant under local diffeomorphisms (see [PIZ13, §2.22]), it must be constant on each Klein stratum. □

## II. DIFFEOLOGICAL STRUCTURE OF ORBIT SPACES AND DIMENSION

In this section, we clarify the local structure of the diffeological quotient of a connected manifold M by the smooth action of a compact Lie group G. We then introduce the (intrinsic) diffeological dimension map on the orbit space M/G and establish a fundamental formula relating it to the dimensions of the manifold and its orbits.

**5. Local Structure of Orbit Spaces.**

Let G be a compact Lie group acting smoothly on a connected manifold M. Let $\mathcal{O} = G \cdot x$ be the orbit of a point $x \in M$, and let $H = \mathrm{Stab}_G(x)$ be its stabilizer.

The classical Slice Theorem (see, e.g., [Pal60, Bre72]) states that there exists a G-invariant neighbourhood $V \subset M$ of the orbit $\mathcal{O}$, and a vector space E (the normal slice at $x$) on which H acts orthogonally, such that V is G-equivariantly diffeomorphic to the associated bundle $G \times_H E = (G \times E)/H$.

When we equip M/G with the quotient diffeology, this local structure translates directly into a local diffeomorphic equivalence:

**Proposition (Diffeological Local Model).** *As a diffeological space, the quotient space* M/G *is locally diffeomorphic to some quotient* E/H *where* E *is an Euclidean space and* H *is a compact subgroup of* O(E) *acting orthogonally. Precisely, if* $V \subset M$ *is a* G*-invariant open neighborhood of an orbit* $\mathcal{O} = G \cdot x$ *and* V *is* G*-equivariantly diffeomorphic to* $G \times_H E$ *(where* $H = \mathrm{Stab}_G(x)$*), then the quotient space* V/G *(with the subspace diffeology from* M/G*) is diffeomorphic to* E/H *(with the quotient diffeology).*

*Proof.* The G-equivariant diffeomorphism $\Phi : V \to G \times_H E$ induces a diffeomorphism $\bar{\Phi} : V/G \to (G \times_H E)/G$. Thus, we need to show that $(G \times_H E)/G$ is diffeomorphic to E/H. By definition, $G \times_H E = (G \times E)/H$. So, $(G \times_H E)/G = ((G \times E)/H)/G$. The group G acts on $G \times E$ by $g' \cdot (g, \xi) = (g'g, \xi)$. The group H acts on $G \times E$ by $h \cdot (g, \xi) = (gh^{-1}, h\xi)$. These two actions commute: $g' \cdot (h \cdot (g, \xi)) = g' \cdot (gh^{-1}, h\xi) = (g'(gh^{-1}), h\xi) = ((g'g)h^{-1}, h\xi)$. $h \cdot (g' \cdot (g, \xi)) = h \cdot (g'g, \xi) = ((g'g)h^{-1}, h\xi)$. Since the actions of G and H on $G \times E$ commute, the iterated quotients are diffeomorphic: $((G \times E)/H)/G \simeq ((G \times E)/G)/H$. The quotient $(G \times E)/G$ is diffeomorphic to E via the projection $\mathrm{pr}_E : G \times E \to E$, $\mathrm{pr}_E(g, \xi) = \xi$. This map is G-invariant in the first factor and induces a diffeomorphism $(G \times E)/G \to E$. Therefore, $((G \times E)/G)/H \simeq E/H$. Combining the diffeomorphisms, we have $V/G \simeq (G \times_H E)/G \simeq E/H$. □

**6. The Orthofold Category.**

The local model E/H, where E is a Euclidean space and $H \subset O(E)$ is a compact subgroup acting orthogonally, appears fundamental for orbit spaces of compact group actions on manifolds. This suggests the definition of a class of diffeological spaces, tentatively termed *orthofolds*, characterized by being locally diffeomorphic, at every point, to such a quotient E/H, in the spirit of [PIZ13, Section "Modeling Diffeology"]. Formally:

**Definition (Orthofolds).** *A diffeological space* X *is called an* orthofold *if* X *is locally diffeomorphic, at each point* $x$*, to some quotient* E/H*, where* E *is an Euclidean vector space and* H *a subgroup of the orthogonal group* O(E)*. The pair* (E, H) *is called a* local model *for the orthofold at* $x$.

This category would naturally include smooth manifolds (where H is trivial) and orbifolds (where H is finite). All orbit spaces M/G of manifolds by compact Lie group actions are natural examples of orthofolds.

**Proposition (Orbifolds within Orthofolds).** *A diffeological space* X *is an orthofold with a constant dimension map if and only if it is an orbifold. So, an orbifold is an orthofold of constant dimension.*

*Proof.* Let X be an orthofold. By definition, for any $x \in X$, there is a local model diffeomorphic to E/H where E is a Euclidean space and $H \subset O(E)$ is a compact subgroup. The diffeological dimension of this local model at the point corresponding to $x$ (the image of $0_E$) is dim(E) (by the Lemma in the proof of dimension formula).

First: If X is an $n$-orbifold, it is, by definition (e.g., as a diffeological space [IKZ10]), locally diffeomorphic to E/H where dim(E) = $n$ and H is a *finite* subgroup of O(E). The diffeological dimension of such a local model E/H (with H finite) is dim(E) = $n$ at every point within that model (since finite group actions on E have zero-dimensional orbits outside the origin, and the dimension formula $\dim_{[v]_H}(E/H) = \dim(E) - \dim(H \cdot v)$ applies). Since this holds for all local models, an $n$-orbifold has a constant dimension map equal to $n$.

Next: Assume X is an orthofold with a constant dimension map, say $\dim_y(X) = d$ for all $y \in X$. Consider any local model E/H for X. We must have dim(E) = $d$. For any point $[v]_H \in E/H$ (where $v \in E$), the dimension of E/H at this point is given by the formula $\dim_{[v]_H}(E/H) = \dim(E) - \dim(\mathcal{O}_v)$, where $\mathcal{O}_v = H \cdot v$ is the orbit of $v$ under H within E (this formula is analogous to Proposition "Quotient Dimension-Map Formula" of article 7 applied to the action of H on E). Since the dimension map is constant and equal to $d = \dim(E)$, we must have $\dim(E) - \dim(\mathcal{O}_v) = \dim(E)$ for all $v \in E$. This implies $\dim(\mathcal{O}_v) = 0$ for all $v \in E$. For a compact Lie group H acting orthogonally on E, if all its orbits $H \cdot v$ are zero-dimensional, then H must be a discrete group [PIZ13, §1.81]. Since H is also compact, it must be a finite group. Therefore, every local model for X is of the form E/H with H finite. This means X is an orbifold. □

**Note.** It is worth noting that other classes of diffeological spaces also share the property of having a constant diffeological dimension, such as *quasifolds* [IZP21]. By definition, a quasifold is locally diffeomorphic to a quotient $\mathbf{R}^n/\Gamma$, where $\Gamma$ is a countable subgroup of the affine group Aff($\mathbf{R}^n$) such that the orbits are discrete. At a point $x$ in such a local model, the isotropy group $\Gamma_x = \{\gamma \in \Gamma \mid \gamma(x) = x\}$ is a countable subgroup of Aff($\mathbf{R}^n$) fixing $x$. The action of $\Gamma_x$ on the tangent space at $x$ (which serves as the slice model) is linear, given by the linear part of $\Gamma_x$, a countable subgroup of GL($n, \mathbf{R}$). For a quasifold to be an orthofold, its local models must be diffeomorphic to E/H where E is a Euclidean space and H is a compact subgroup of O(E). However, a countable subgroup of GL($n, \mathbf{R}$) can only be compact if it is finite. Therefore, quasifolds with infinite isotropy groups (i.e., strict quasifolds) cannot be locally diffeomorphic to quotients by compact groups, and thus strict quasifolds are not orthofolds.

## 7. The Diffeological Dimension Map and Formula.

The diffeological dimension map $\dim_x(X)$ is a fundamental intrinsic invariant defined point-wise for any diffeological space X [PIZ07, PIZ13]. For the orbit space M/G, equipped with the quotient diffeology, this map $\mathcal{O} \mapsto \dim_\mathcal{O}(M/G)$ provides diffeologically invariant labels attached to each orbit $\mathcal{O} \in M/G$.

As the diffeological dimension is an invariant under local diffeomorphisms, it is constant on the Klein strata of M/G (i.e., on the orbits of the pseudogroup of local diffeomorphisms of M/G).

The diffeological dimension in M/G captures geometric characteristics of the singularities, reflecting how orbits deviate from being principal. The dimension of M/G at a point $\mathcal{O}$ corresponding to an orbit $G \cdot x$ is given by the following formula:

**Proposition (Quotient Dimension-Map Formula).** *Let* orbit : $M \to M/G$ *be the projection from* M *onto the space of orbits* M/G*, equipped with the quotient diffeology. For a point* $x \in M$*, let* $\mathcal{O} = \text{orbit}(x)$ *be the orbit* $G \cdot x$*. The diffeological dimension of* M/G *at the point* $\mathcal{O}$ *obeys the formula:*

$$\dim_{\mathcal{O}}(M/G) = \dim(M) - \dim(\mathcal{O}). \tag{$\clubsuit$}$$

*Proof.* The formula $\dim_{\mathcal{O}}(M/G) = \dim(M) - \dim(\mathcal{O})$ relies on the local structure of M/G and the dimension of the local models.

**Lemma.** Let E be a finite-dimensional real vector space, and let $H \subset O(E)$ be a compact subgroup of its orthogonal group. Then $\dim_{[0_E]_H}(E/H) = \dim(E)$.

*Proof of the Lemma.* Let $\Delta = E/O(E)$ and $\pi_{O(E)} : E \to \Delta$ be the projection. Let $X = E/H$ and $\pi_H : E \to X$ be the projection. Define $\text{pr} : X \to \Delta$ by $\text{pr}([v]_H) = [v]_{O(E)}$. All three projections are subductions: $\pi_{O(E)}$ and $\pi_H$ by definition of quotient diffeology, and pr because it is induced by the identity on E and respects the equivalence relations. These maps send the origin of E to the respective origins of the quotient spaces: $\pi_H(0_E) = [0_E]_H$, $\pi_{O(E)}(0_E) = [0_E]_{O(E)}$, and $\text{pr}([0_E]_H) = [0_E]_{O(E)}$. Since $\pi_H$ and pr are subductions preserving the origins (in the sense that the origin of the domain maps to the point we are considering as the "origin" in the codomain), by properties of the diffeological dimension under subductions (see [PIZ13, §2.24, 2.26]), we have:

$$\dim_{[0_E]_H}(X) \leq \dim_{0_E}(E), \text{ and } \dim_{[0_E]_{O(E)}}(\Delta) \leq \dim_{[0_E]_H}(X).$$

This yields the chain:

$$\dim_{[0_E]_{O(E)}}(\Delta) \leq \dim_{[0_E]_H}(X) \leq \dim_{0_E}(E).$$

The diffeological dimension of E at any point is its manifold dimension, $\dim_{0_E}(E) = \dim(E)$. From [PIZ13, Exercise 50] (with $E = \mathbf{R}^n$), we know that $\dim_{[0_E]_{O(E)}}(\Delta) = \dim(E)$. Thus,

$$\dim(E) \leq \dim_{[0_E]_H}(X) \leq \dim(E).$$

This forces $\dim_{[0_E]_H}(X) = \dim(E)$, that is, $\dim_{[0_E]_H}(E/H) = \dim(E)$. □

Now, let $\mathcal{O} = \text{orbit}(x) \in M/G$ be the orbit $G \cdot x$, with stabilizer $H = \text{Stab}_G(x)$. By Proposition "Diffeological Local Model" of article 5, M/G is locally diffeomorphic around the point $\mathcal{O} \in M/G$ to E/H, where E is the slice at $x$. The point $\mathcal{O} \in M/G$ corresponds to $[0_E]_H$ in this local model. Since the diffeological dimension is a local invariant, $\dim_{\mathcal{O}}(M/G) = \dim_{[0_E]_H}(E/H)$. By the Lemma, $\dim_{[0_E]_H}(E/H) = \dim(E)$. The dimension of the manifold M at $x$ is related to the dimensions of the orbit and the slice by $\dim(M) = \dim(G \cdot x) + \dim(E)$. The dimension of the orbit $G \cdot x$ is $\dim(G) - \dim(\text{Stab}_G(x)) = \dim(G) - \dim(H)$, but more relevant here is its dimension as a submanifold of M, which is $\dim(\mathcal{O})$. Thus, $\dim(M) = \dim(\mathcal{O}) + \dim(E)$. Substituting $\dim(E) = \dim_{\mathcal{O}}(M/G)$, we get $\dim(M) = \dim(\mathcal{O}) + \dim_{\mathcal{O}}(M/G)$, which yields the formula: $\dim_{\mathcal{O}}(M/G) = \dim(M) - \dim(\mathcal{O})$. □

## 8. Consequences of the Dimension Formula.

The dimension formula (♣) provides insight into the structure of M/G:

**(a)** The set of principal orbits $(M/G)_{pr} = \mathrm{orbit}(M_{pr})$ is the *principal Klein stratum* of M/G. It is an open dense subset of M/G (in the D-topology, which coincides with the standard quotient topology for M/G) and consists of manifold points. All principal orbits $\mathcal{O}_{pr}$ have the same dimension. Thus, the dimension map is constant on $(M/G)_{pr}$, with value $d_{pr} = \dim(M) - \dim(\mathcal{O}_{pr})$, which is consistent with the general property in diffeology that the dimension map is constant on Klein strata.

**(b)** For any singular orbit $\mathcal{O}_s$, its dimension $\dim(\mathcal{O}_s)$ is less than or equal to the dimension of a principal orbit $\dim(\mathcal{O}_{pr})$. This is because the stabilizer of a singular point contains a conjugate of a principal stabilizer, so $\dim(\mathrm{Stab}_G(s)) \geq \dim(\mathrm{Stab}_G(x_{pr}))$, which implies $\dim(\mathcal{O}_s) = \dim(G) - \dim(\mathrm{Stab}_G(s)) \leq \dim(G) - \dim(\mathrm{Stab}_G(x_{pr})) = \dim(\mathcal{O}_{pr})$. By formula (♣), $\dim_{\mathcal{O}_s}(M/G) = \dim(M) - \dim(\mathcal{O}_s) \geq \dim(M) - \dim(\mathcal{O}_{pr}) = d_{pr}$. Thus, $d_{pr}$ is the minimal dimension in M/G.

**(c)** Equality, $\dim_{\mathcal{O}_s}(M/G) = d_{pr}$, occurs for singular orbits $\mathcal{O}_s$ if and only if $\dim(\mathcal{O}_s) = \dim(\mathcal{O}_{pr})$. Such orbits are often referred to as *exceptional orbits*. This happens when the stabilizer $\mathrm{Stab}_G(s)$ has the same dimension as a principal stabilizer but is not conjugate to it.

## 9. Application to Toric Actions and Depth.

The diffeological dimension map provides a natural way to quantify local complexity in quotient spaces arising from group actions, sometimes aligning with or refining existing topological or geometric invariants. Consider an effective action of the $n$-torus $T^n$ on a $2n$-dimensional symplectic manifold $(M^{2n}, \omega)$. Such actions are central to symplectic geometry, and their orbit spaces $Q_n = M^{2n}/T^n$ have a well-studied combinatorial and topological structure.

In many local models for these quotient spaces (e.g., near a fixed point or a point with a non-trivial stabilizer), the structure is locally diffeomorphic to a quotient of the form E/H where E is a Euclidean space and H is a compact subgroup of O(E) acting orthogonally. A standard example of such a local model is $\mathbf{C}^n/T^n$, where $T^n$ acts by $(e^{i\theta_1}, \ldots, e^{i\theta_n}) \cdot (z_1, \ldots, z_n) = (e^{i\theta_1} z_1, \ldots, e^{i\theta_n} z_n)$. The quotient $\mathbf{C}^n/T^n$ is diffeomorphic to the positive orthant $[0, \infty[^n$ via the map $(z_1, \ldots, z_n) \mapsto (|z_1|^2, \ldots, |z_n|^2)$.

In the topological analysis of such spaces, a label called the *depth* is often associated with points in the quotient, see for example [KS]. For a point $t = (t_1, \ldots, t_n) \in [0, \infty[^n$ representing an orbit, its depth is typically defined as the number of coordinates $t_i$ that are equal to zero. This corresponds to the dimension of the subtorus in $T^n$ that stabilizes points in $\mathbf{C}^n$ projecting to $t$.

Let $t = (|z_1|^2, \ldots, |z_n|^2)$ be a point in $\mathbf{C}^n/T^n$. According to the dimension formula (♣), the diffeological dimension at $t$ is $\dim_t(\mathbf{C}^n/T^n) = \dim(\mathbf{C}^n) - \dim(\mathcal{O}_z)$, where $\mathcal{O}_z = T^n \cdot z$. The dimension of the orbit $\mathcal{O}_z$ is $\dim(T^n) - \dim(\mathrm{Stab}_{T^n}(z))$. The stabilizer $\mathrm{Stab}_{T^n}(z)$ is isomorphic to $T^k$ if exactly $k$ coordinates of $z$ are zero. Thus, $\dim(\mathrm{Stab}_{T^n}(z)) = k$, and $\dim(\mathcal{O}_z) = n - k$. The integer $k$ is precisely the depth of the point $t$.

Substituting this into the formula:

$$\dim_t(\mathbf{C}^n/\mathrm{T}^n) = \dim(\mathbf{C}^n) - (n-k) = 2n - (n-k) = n+k.$$

Thus, the depth, an important combinatorial invariant in the study of toric symplectic manifolds, is directly related to the diffeological dimension of this local model of the orbit space by the formula:

$$\dim_t(\mathbf{C}^n/\mathrm{T}^n) = n + \mathrm{depth}(t).$$

This shows how the diffeological dimension naturally incorporates and quantifies the "degree of singularity" captured by the depth in this context. For a point $t$ in the interior of $[0,\infty[^n$ (depth $k=0$), $\dim_t(\mathbf{C}^n/\mathrm{T}^n) = n$. For a point at the origin of $[0,\infty[^n$ (depth $k=n$), $\dim_t(\mathbf{C}^n/\mathrm{T}^n) = 2n = \dim(\mathbf{C}^n)$.

A final remark: It is important to distinguish this quotient diffeology on $[0,\infty[^n$ (obtained from $\mathbf{C}^n/\mathrm{T}^n$) from its standard diffeology as a manifold with corners when viewed as a subset of $\mathbf{R}^n$. For the latter, the diffeological dimension is $\infty$ at any point on the boundary, and $n$ for interior points. In contrast, the quotient diffeology yields $\dim_t(\mathbf{C}^n/\mathrm{T}^n) = n + \mathrm{depth}(t)$ for all $t$, reflecting the structure inherited from $\mathbf{C}^n$ and the $\mathrm{T}^n$-action. Adding the depth of the strata to the topological structure of a manifold with corners for these quotients mixes concepts from two different categories. It is certainly more satisfactory for the mind to stay within a single category, and considering the quotient as a diffeological space provides a perfectly coherent framework.

## III. Examples

This section presents several examples of quotients of manifolds by compact group actions. They are selected to illustrate the intricate, and often surprising, relationship between the classical orbit-type stratification on M and the intrinsic Klein stratification of the quotient space M/G. Each example reveals a subtlety that challenges the simple expectation of a direct correspondence and motivates the more refined analysis presented in the subsequent sections.

### 10. Example $(\mathrm{S}^2 \times \mathrm{S}^2)/\mathrm{SO}(3)$.

This first example illustrates a baseline case where the relationship appears simple. The orbit-type stratification of M projects neatly onto the Klein stratification of the quotient.

We consider the group SO(3) acting diagonally on $M = \mathrm{S}^2 \times \mathrm{S}^2$. The stabilizers are:

$$\mathrm{Stab}_{\mathrm{SO}(3)}(x,y) = \begin{cases} \{\mathrm{id}\} & \text{if } x \neq y \text{ and } x \neq -y \text{ (principal type)}, \\ \mathrm{SO}(2,x) & \text{if } x = y \text{ or } x = -y \text{ (singular type)}, \end{cases}$$

where $\mathrm{SO}(2,x)$ is the subgroup fixing $x$. All singular stabilizers are conjugate to SO(2). Thus, the orbit-type stratification of M consists of two strata: a principal stratum $\mathrm{M}_{(\{\mathrm{id}\})}$ and a singular stratum $\mathrm{M}_{(\mathrm{SO}(2))}$.

The orbit space $Q = M/\mathrm{SO}(3)$ is diffeomorphic to the interval $[-1,+1]_\pi$, where the map $\pi(x,y) = \langle x, y \rangle$ distinguishes the orbits. The Klein stratification of Q consists of two strata: the interior $K_1 = ]-1,1[$ and the boundary $K_2 = \{\pm 1\}$.

The projection of the orbit-type strata aligns perfectly: $\mathrm{orbit}(M_{(\{\mathrm{id}\})}) = K_1$ and $\mathrm{orbit}(M_{(\mathrm{SO}(2))}) = K_2$. This example suggests a simple correspondence might hold in general.

## 11. Example $P^2(\mathbf{R})/\mathrm{SO}(2)$.

This example shows that different types of singularities can coexist within a simple quotient space. Consider the action of $G = \mathrm{SO}(2, k)$ (rotations around the unit vector $k$) on the real projective plane $M = P^2(\mathbf{R}) = S^2/\{\pm 1\}$. The orbit space $Q = M/G$ is diffeomorphic to the interval $[0,1]_\pi$, with $\pi : [x] \mapsto \langle x, k \rangle^2$ realizing the quotient by pushforward.

An analysis of the stabilizers and slice representations reveals three distinct local models in the quotient:

* The interior $]0,1[$ consists of principal points, with local model $\mathbf{R}$.
* The endpoint $t = 0$ corresponds to an exceptional orbit type. Its local model is $\mathbf{R}/\mathrm{O}(1) = \Delta_1$, making it an **orbifold point**. Its diffeological dimension is 1.
* The endpoint $t = 1$ corresponds to a fixed point. Its local model is $\mathbf{R}^2/\mathrm{O}(2) = \Delta_2$, making it a more general **orthofold point**. Its diffeological dimension is 2.

In this case, the three orbit-type strata project to three distinct Klein strata ($\{0\}$, $\{1\}$, and $]0,1[$). The key pedagogical feature is that the quotient space is a simple segment whose endpoints represent fundamentally different kinds of singularities, a distinction clearly captured by the diffeology.

## 12. Example $P^2(\mathbf{C})/\mathrm{SO}(3)$.

This example immediately complicates the picture, showing that distinct orbit-type strata in M can map into the same Klein stratum in the quotient.

Consider the complex projective plane $M = P^2(\mathbf{C})$ with a specific action of $G = \mathrm{SO}(3)$ (see [PI91]). The stabilizer conjugacy classes that occur define three orbit-type strata:

* A principal stratum $M_{(\mathbf{Z}_2)}$.
* A singular stratum $M_{(\mathrm{SO}(2))}$.
* Another singular stratum $M_{(\mathrm{O}(2))}$.

The groups SO(2) and O(2) are not conjugate in SO(3), so these are genuinely different orbit types. The orbit space $Q = M/\mathrm{SO}(3)$ is again diffeomorphic to an interval, say $[0,1]_\pi$, with $\pi : [X, Y] \mapsto 1 - 4\|X \wedge Y\|^2$ the projection from $P^2(\mathbf{C})$ onto $[0, 1]$, realizing the quotient by pushforward.

The orbit-type strata project as follows: $\mathrm{orbit}(M_{(\mathbf{Z}_2)}) = ]0,1[$, $\mathrm{orbit}(M_{(\mathrm{SO}(2))}) = \{0\}$, and $\mathrm{orbit}(M_{(\mathrm{O}(2))}) = \{1\}$.

However, an analysis of the local diffeological structure reveals a surprise. The local models are:

* For $t \in ]0,1[$: model $\mathbf{R}$.
* For $t = 0$: model is diffeomorphic to $\mathbf{R}^2/\mathrm{O}(2) = \Delta_2$.
* For $t = 1$: model is also $\mathbf{R}^2/\mathrm{O}(2) = \Delta_2$.

Since the local models at $t = 0$ and $t = 1$ are diffeomorphic, these two points belong to the *same* Klein stratum. The Klein stratification of Q consists of only two strata: $K_1 = ]0, 1[$ and $K_2 = \{0, 1\}$.

This presents a puzzle for the orbit-type correspondence: two distinct singular orbit-type strata, $M_{(SO(2))}$ and $M_{(O(2))}$, both project into the singular Klein stratum $K_2$. The orbit-type stratification on M is finer than what is needed to describe the intrinsic geometry of M/G.

## 13. Example: $P^2(\mathbf{C})/U(1)$.

This final example reveals an even deeper subtlety: a single orbit-type stratum in M can map to points in *different* Klein strata in the quotient. This demonstrates conclusively that the orbit-type stratification is insufficient.

This example is from [Aud91]. Consider the action of $G = U(1)$ on $M = P^2(\mathbf{C})$ given by $\tau \cdot [z_1 : z_2 : z_3] = [z_1 : \tau z_2 : \tau^2 z_3]$. The stabilizers that occur are $\{1\}$, $\{\pm 1\}$, and $U(1)$. This gives three orbit-type strata: $M_{(\{1\})}$, $M_{(\{\pm 1\})}$, and the fixed-point stratum $M_{(U(1))}$, which consists of the three points $P_0 = [1 : 0 : 0]$, $P_1 = [0 : 1 : 0]$, and $P_2 = [0 : 0 : 1]$.

All three fixed points belong to the same orbit-type stratum, $M_{(U(1))}$. One might expect their images in the quotient, orbit($P_0$), orbit($P_1$), and orbit($P_2$), to belong to the same Klein stratum. However, the local structure is determined by the slice representation at each point:

* At $P_0$: Slice action $\rho_0(\tau) \cdot (w_1, w_2) = (\tau w_1, \tau^2 w_2)$.
* At $P_1$: Slice action $\rho_1(\tau) \cdot (w_1, w_2) = (\bar{\tau} w_1, \tau w_2)$.
* At $P_2$: Slice action $\rho_2(\tau) \cdot (w_1, w_2) = (\bar{\tau}^2 w_1, \bar{\tau} w_2)$.

The representations $\rho_0$ and $\rho_2$ are equivalent, but $\rho_1$ is not equivalent to them. Consequently, the local models for orbit($P_0$) and orbit($P_2$) are diffeomorphic to each other, but this will not be the case for the local model for orbit($P_1$). Indeed, the slice actions at $P_0$ and $P_2$ admit non-trivial finite stabilizers on subspaces of the slice, which implies that in the quotient, the points orbit($P_0$) and orbit($P_2$) lie in the closure of the exceptional (orbifold) stratum orbit($M_{\{\pm 1\}}$). Conversely, the slice action at $P_1$ is free away from the origin, meaning orbit($P_1$) is topologically isolated from the exceptional stratum. Therefore, orbit($P_0$) and orbit($P_2$) lie in one Klein stratum, while orbit($P_1$) lies in a completely different one. Here, a single orbit-type stratum (the discrete set $M_{(U(1))}$) is "torn apart" by the projection to the quotient space.

These examples collectively show that to understand the Klein stratification of M/G, we need a partition of M that is finer than the orbit-type stratification and sensitive to the specific stabilizer and its slice representation, not just its conjugacy class. This motivates the introduction of the isostabilizer decomposition, which we will define and analyze in the following sections.

## 14. Example: Orbifolds.

The case of orbifolds provides a useful classical benchmark. An $n$-orbifold is a diffeological space locally diffeomorphic to $\mathbf{R}^n/\Gamma$ for some finite subgroup $\Gamma \subset GL(n, \mathbf{R})$. The conjugacy class of $\Gamma$ is the *structure group* at that point.

A key result from [IKZ10] states that two points in an orbifold belong to the same Klein stratum if and only if their structure groups are conjugate.

When an orbifold arises as a quotient M/G of a manifold by a finite group G, the structure group at orbit($x$) is precisely the stabilizer $\mathrm{Stab}_G(x)$. In this situation, the Klein stratification of the quotient M/G corresponds exactly to the projection of the orbit-type stratification of M. This aligns with the classical understanding of orbifolds. The direct correspondence between the orbit-type stratification of M and the Klein stratification of M/G holds in this setting. This serves as a crucial benchmark, highlighting that the more complex behaviors observed in our other examples—where the correspondence breaks down—are characteristic of actions involving stabilizers with non-trivial identity components (such as SO(2) and U(1)).

## IV. THE ISOSTABILIZER DECOMPOSITION

The examples in the previous section demonstrate that the classical orbit-type stratification is not always the correct tool to predict the local structure of the quotient space M/G. To establish a precise relationship, we must introduce a finer partition of M that is sensitive not only to the stabilizer subgroup itself but also to the equivalence class of its action on the normal slice. This section formally defines this partition.

Let G be a compact Lie group acting smoothly on a connected smooth manifold M.

### 15. Orbits and Stabilizers.

For any point $x \in M$, its *orbit* is the set $G \cdot x = \{g \cdot x \mid g \in G\}$. The *stabilizer* (or isotropy group) of $x$ is the subgroup $\mathrm{Stab}_G(x) = \{g \in G \mid g \cdot x = x\}$. The stabilizer is a closed subgroup of G, and all points on the same orbit have conjugate stabilizers: $\mathrm{Stab}_G(g \cdot x) = g\,\mathrm{Stab}_G(x)g^{-1}$.

### 16. Isostabilizer Submanifolds.

For each closed subgroup $K \subseteq G$, we define the *isostabilizer set* $M_K = \{x \in M \mid \mathrm{Stab}_G(x) = K\}$. The collection of non-empty sets $\{M_K\}$ for all closed subgroups $K \subseteq G$ forms a partition of M. Each of these sets is, in fact, a submanifold.

**Theorem (Isostabilizer Sets are Submanifolds).** *For any closed subgroup* $K \subseteq G$, *the set* $M_K$ *is a smooth submanifold of* M *(possibly empty).*

*Proof.* By the Differentiable Slice Theorem [Pal60, Bre72], for any $x \in M_K$, a G-invariant open neighborhood U of the orbit $G \cdot x$ is G-equivariantly diffeomorphic to an associated bundle $G \times_K B$, where B is a K-invariant open ball in the normal space $N_x$ on which K acts linearly. The subset $M_K \cap U$ corresponds under this diffeomorphism to the points $[g, v]_K \in G \times_K B$ with stabilizer K. This condition holds if and only if $g \in N(K)$ (the normalizer of K) and $v \in B^K$ (the set of K-fixed points in B).

Thus, $M_K \cap U$ is diffeomorphic to the bundle $N(K) \times_K B^K$, which is trivial and hence diffeomorphic to the product $(N(K)/K) \times B^K$. The quotient $N(K)/K$ is a manifold. Since the action of K on $N_x$ is linear, its fixed-point set $(N_x)^K$ is a linear subspace. The set $B^K = B \cap (N_x)^K$ is an open subset of this subspace and is therefore a smooth manifold.

The product of two manifolds is a manifold, proving that $M_K$ is locally a manifold, and thus is a manifold itself. □

**Remark.** While this is a standard result, we include the proof to highlight the local product structure $(N(K)/K) \times B^K$. This structure is the foundational brick for the fine analysis of the quotient diffeology, which is not an obvious consequence of the coarser orbit-type stratification.

### 17. The Isostabilizer Decomposition.

The local structure of the quotient space $M/G$ at a point orbit($x$) is determined by the action of the stabilizer $\text{Stab}_G(x)$ on the slice at $x$. As shown in Example 13, points with the same stabilizer K can still produce different local structures in the quotient if their slice representations are not equivalent. This motivates looking at the sets where the slice representation remains constant.

**Proposition (Slice Equivalence on Components).** *Let* C *be a connected component of an isostabilizer submanifold* $M_K$. *For any two points* $x, x' \in C$, *the slice representations* $\rho_x : K \to O(N_x)$ *and* $\rho_{x'} : K \to O(N_{x'})$ *are equivalent as representations of* K.

*Proof.* Let C be a connected component of $M_K$. We can equip M with a G-invariant Riemannian metric. The normal spaces $N_x$ to the orbits $G \cdot x$ for $x \in C$ form a smooth vector bundle $\mathcal{N} \to C$. The slice representation at $x$ is the action of K on the fiber $N_x$. Since C is connected, we can connect any two points $x, x'$ by a path. The continuous dependence of the representation along this path, combined with the fact that the space of equivalence classes of representations of a compact group on a fixed-dimensional vector space is discrete, implies that the representations $\rho_x$ and $\rho_{x'}$ must be equivalent (see, e.g., [BT85, Fol95]). □

This proposition shows that the connected components of the isostabilizer submanifolds are the fundamental "atoms" of the action, on which the local data (stabilizer and slice representation) is constant. This justifies the following central definition.

**Definition (Isostabilizer Decomposition).** *The* isostabilizer decomposition *of the* G-*manifold* M *is the partition of* M *into the connected components of the isostabilizer submanifolds* $\{M_K\}$. *We denote the set of these components by* $\text{Iso}(M, G)$.

This decomposition is finer than the classical orbit-type stratification. As we will show in the next section, it is precisely this decomposition of M that maps coherently to the Klein stratification of the quotient space $M/G$.

## V. THE MAIN CORRESPONDENCE THEOREM

The previous sections have established the necessary tools: the intrinsic Klein stratification of the quotient space $M/G$ and the refined isostabilizer decomposition of the manifold M. We are now in a position to state and prove the precise relationship between them. We will show that the isostabilizer decomposition of M is exactly the structure that maps coherently to the Klein stratification of $M/G$.

### 18. The Induced Map on Decompositions.

Recall that Iso(M, G) is the set of connected components of the isostabilizer submanifolds $M_K$, and $K\ell(M/G)$ is the set of Klein strata of the quotient space. By Proposition "Slice Equivalence on Components" of article 17, for any component $C \in \text{Iso}(M, G)$, the stabilizer K is constant and the slice representation is fixed up to equivalence. This implies that the local model for the quotient, $E_x/K$, is of a constant diffeomorphic type for all points $x \in C$.

Since the Klein strata are, by definition, the sets of points with diffeomorphic local models, the image of any component C under the orbit projection must lie entirely within a single Klein stratum. This observation leads to our main theorem.

**Theorem (The Correspondence Theorem).** *The canonical projection* $\text{orbit} : M \to M/G$ *induces a well-defined, surjective map*

$$\text{orbit}_* : \text{Iso}(M, G) \to K\ell(M/G)$$

*from the isostabilizer decomposition of* M *to the Klein stratification of* M/G.

*Proof.* **Well-definedness:** Let $C \in \text{Iso}(M, G)$. By Proposition "Slice Equivalence on Components", for any two points $x_1, x_2 \in C$, their stabilizers are the same group K, and their slice representations are equivalent. The local model of the quotient M/G at a point orbit($x$) is determined by the equivalence class of the slice representation of $\text{Stab}_G(x)$. Since this is constant for all $x \in C$, the points orbit($x_1$) and orbit($x_2$) have diffeomorphic local models. Therefore, they belong to the same Klein stratum. This shows that the entire set orbit(C) is contained within a single Klein stratum, so the map $\text{orbit}_*$ is well-defined.

**Surjectivity:** Let $\mathcal{K} \in K\ell(M/G)$ be an arbitrary Klein stratum. By definition, there must exist at least one point $\mathcal{O} \in \mathcal{K}$. This point is an orbit, so $\mathcal{O} = \text{orbit}(x)$ for some $x \in M$. This point $x$ must belong to some isostabilizer submanifold $M_K$, and therefore to one of its connected components, say $C \in \text{Iso}(M, G)$. By construction, $\text{orbit}_*(C)$ is the Klein stratum containing orbit($x$), which is $\mathcal{K}$. Thus, every Klein stratum is in the image of $\text{orbit}_*$. □

### 19. The Inverse Klein Stratification.

The Correspondence Theorem allows us to define a new, canonical stratification on M induced conversely by the quotient's geometry.

Let $\{K\ell_i\}$ be the Klein stratification of M/G. We define the **Inverse Klein Stratification** of M to be the partition $\{S_i\}$ where each set $S_i = \text{orbit}^{-1}(K\ell_i)$ is the preimage of a Klein stratum.

**Proposition (The Inverse Klein is a Stratification).** *The Inverse Klein Stratification* $\{S_i\}$ *is a true stratification of the manifold* M.

*Proof.* Let $K\ell_j \subset \overline{K\ell_i}$ in M/G. We must show that $S_j \subset \overline{S_i}$. Let $y \in S_j$. For any open neighborhood W of $y$, orbit(W) is an open neighborhood of orbit($y$). Since orbit($y$) $\in \overline{K\ell_i}$, orbit(W) must intersect $K\ell_i$. This implies there exists an $x \in W$ such that orbit($x$) $\in K\ell_i$, meaning $x \in S_i$. Thus, any neighborhood of $y$ intersects $S_i$, so $y \in \overline{S_i}$. □

Our key examples show that the relationship between this new stratification $\{S_i\}$ and the classical orbit-type stratification $\{M_{(H)}\}$ is one of **incomparability**.

- In Example 12, the Inverse Klein Stratification is strictly **coarser** than the orbit-type stratification, as it merges the distinct orbit types (SO(2)) and (O(2)).
- In Example 13, the Inverse Klein Stratification is strictly **finer** than the orbit-type stratification, as it splits the single orbit type (U(1)) into two strata.

This clarifies the hierarchy of partitions on M: the isostabilizer decomposition $\{C_\alpha\}$ is the finest of the three, refining both the classical orbit-type stratification and the new Inverse Klein Stratification. This solidifies its role as the fundamental "atomic" partition for understanding the quotient's geometry.

## 20. Characterization of Singularity Types in M/G.

The Correspondence Theorem allows us to classify the geometric nature of the points in the quotient space by examining the stabilizers of the corresponding points in the manifold.

**Proposition (Geometric Nature of Quotient Points).** *Let $\mathscr{O} = \mathrm{orbit}(x) \in M/G$ be a point in the orbit space, with stabilizer $K = \mathrm{Stab}_G(x)$.*

* *If $x$ is a principal point, then $\mathscr{O}$ is a **manifold point**. Its local model is a Euclidean space.*
* *If $x$ is an exceptional singular point (i.e., $\dim(K) = \dim(K_{\mathrm{pr}})$ but $K$ is not conjugate to $K_{\mathrm{pr}}$), then $\mathscr{O}$ is a non-trivial **orbifold point**. Its local model is $E/\Gamma$ where $\Gamma$ is a non-trivial finite group.*
* *If $x$ is a non-exceptional singular point (i.e., $\dim(K) > \dim(K_{\mathrm{pr}})$), then $\mathscr{O}$ is a **general orthofold point**. Its local model is $E/K$ where the effective action of $K$ is by a non-finite compact Lie group.*

*Proof.* The local model at $\mathscr{O}$ is determined by the effective action of $K$ on the slice $E$. *Principal:* The effective action of a principal stabilizer on its slice is trivial [Bre72, IV, 3.2(iii)]. The local model is $E/\{1\} \simeq E$. *Exceptional:* The condition $\dim(K) = \dim(K_{\mathrm{pr}})$ implies the effective action on the slice is by the finite group $\Gamma = K/K^0$. Since the point is singular, this action is non-trivial. The local model is $E/\Gamma$, which is an orbifold. *Non-exceptional Singular:* The condition $\dim(K) > \dim(K_{\mathrm{pr}})$ implies the effective action on the slice has a non-trivial identity component. The local model is therefore an orthofold but not an orbifold. □

## 21. The Meaning of Non-Injectivity.

The Correspondence Theorem establishes a surjective map $\mathrm{orbit}_* : \mathrm{Iso}(M, G) \to K\ell(M/G)$. As our examples in Section III conclusively demonstrate, this map is generally not injective. This non-injectivity is not a flaw in the analysis, but rather a deep insight into the nature of diffeological quotients.

Non-injectivity occurs when distinct components $C_1, C_2 \in \mathrm{Iso}(M, G)$ map to the same Klein stratum. This happens if and only if their corresponding local models are diffeomorphic, even if the components themselves are profoundly different from the perspective of the action on M. We saw two primary ways this can happen:

(1) **Different Stabilizers, Same Local Model:** In Example 12, components of $M_{\mathrm{SO}(2)}$ and $M_{\mathrm{O}(2)}$ mapped to the same Klein stratum. The stabilizers SO(2) and O(2) are not even conjugate, but their actions on their respective slices produce diffeomorphic quotients.
(2) **Same Stabilizer, Different Local Models:** In Example 13, different connected components of the single set $M_{U(1)}$ mapped to *different* Klein strata. This happened because the slice representations at these fixed points were not equivalent: they do not define local diffeomorphic quotient.

This leads to a fundamental conclusion: the intrinsic smooth identity of the quotient space M/G is determined entirely by the set of local models {E/H} that appear. The Klein stratification captures this intrinsic identity. The specific details of the action on M—such as the conjugacy classes of stabilizers or the global structure of the isostabilizer sets—can be "forgotten" in the process of taking the quotient. The quotient space only remembers the local diffeomorphic types of its singularities, the local geometry.

## VI. Concluding Remarks

In this paper, we have clarified the precise relationship between the smooth action of a compact Lie group G on a manifold M and the intrinsic smooth structure of the resulting orbit space M/G. By employing the tools of diffeology, we have moved beyond the classical orbit-type stratification to provide a complete picture of how the local geometry of the action translates into the local geometry of the quotient.

Our investigation yielded three main results. First, we identified the correct partition on the manifold M for this purpose: not the classical orbit-type stratification, but the finer *isostabilizer decomposition*, whose elements are the connected components of submanifolds where the stabilizer is constant. We showed that on each such component, the slice representation is fixed up to equivalence, making these components the fundamental "atoms" of the group action.

Second, we established the *Correspondence Theorem*, which states that the canonical projection induces a surjective map from this isostabilizer decomposition of M to the Klein stratification of M/G. This theorem provides the definitive dictionary between the structure of the action and the intrinsic singularity types of the quotient.

Third, and perhaps most profoundly, we explored the meaning of this map's *non-injectivity*. Our examples demonstrate that distinct structural elements on M can map to diffeologically identical points in the quotient. This reveals a fundamental principle: the quotient space M/G possesses an intrinsic geometric identity that can be independent of the specific details of its construction. It "forgets" the global structure of the action on M, remembering only the local quotient geometry defined by the local diffeomorphisms, which in turn captures the diffeomorphic types of its singularities.

This perspective gives rise to the concept of an *orthofold*, a diffeological space that is locally diffeomorphic to a quotient E/H of a Euclidean space by a compact orthogonal group action. Our analysis shows that orbit spaces of compact group actions are the canonical examples of orthofolds. This class of spaces naturally contains both smooth manifolds (H is trivial) and orbifolds (H is finite) as special cases. Our characterization of the singularity types in Proposition "Geometric Nature of Quotient Points" of article 20 provides a clear hierarchy within this class: principal orbits correspond to manifold points, exceptional orbits to orbifold points, and other singular orbits to general orthofold points. The study of orthofolds as a distinct and natural class of geometric objects, motivated by the structure of orbit spaces, presents a fruitful area for future investigation. Furthermore, the set Iso(M, G) itself can be endowed with a quotient diffeology, turning the map orbit$_*$ into a subduction and opening new avenues for studying the global structure of these decompositions.

The power of the diffeological framework is not limited to such quotients. As illustrated by the space of geodesics on a torus, Geod($T^2$) [PIZ25, §72],[2] diffeology provides essential tools for analyzing the intrinsic smooth structure of a wide range of singular spaces (not limited to quotients of manifolds by compact Lie groups) where traditional differential geometry may lack a suitable global language and tools.

These results underscore the utility of diffeology in providing a consistent and nuanced framework for the study of singular spaces. By capturing the complete smooth structure that remains after reduction, the diffeological approach offers a richer understanding of a quotient's geometry than its topology alone. The study of orthofolds as a distinct and natural class of geometric objects, motivated by the structure of orbit spaces, presents a fruitful area for future investigation.

---

[2]Quotient of a manifold ($T^3$) by a 1-dimensional foliation with closed and non-closed leaves.

SERAP GÜRER, GALATASARAY UNIVERSITY, ORTAKÖY, ÇIRAĞAN CD. NO:36, 34349 BEŞIKTAŞ / İSTANBUL, TURKEY.

*Email address*: `sgurer@gsu.edu.tr`

PATRICK IGLESIAS-ZEMMOUR, EINSTEIN INSTITUTE OF MATHEMATICS, THE HEBREW UNIVERSITY OF JERUSALEM, EDMOND J. SAFRA CAMPUS, GIVAT RAM, 9190401 JERUSALEM, ISRAEL.

*Email address*: `piz@math.huji.ac.il`

*URL*: `http://math.huji.ac.il/~piz`